\documentclass[a4paper,11pt,draft]{amsart}

\usepackage{amsmath}
\usepackage{amsthm}
\usepackage{amsfonts}
\usepackage{amssymb}
\usepackage[all]{xy}

\newcommand{\ZZ}{\mathbb{Z}}

\newtheorem{theorem}{Theorem}
\newtheorem{proposition}{Proposition}

\newtheorem{definition}{Definition}[theorem]
\newtheorem{lemma}{Lemma}[theorem]

\newtheorem{example}{Example}

\begin{document}

\title{Quiver Grassmannians associated with  string modules}

\author{Giovanni Cerulli Irelli}
\address{Universit\`a degli studi di Padova, Dipartimento di Matematica Pura ed Applicata. Via Trieste 63, 35121, Padova (ITALY)}
\email{giovanni.cerulliirelli@gmail.com}


\thanks{Research supported by grant CPDA071244/07 of Padova University}

\begin{abstract}
We provide a technique to compute the Euler--Poincar\'e characteristic of  a  class of projective varieties  called quiver Grassmannians.   This technique applies to quiver Grassmannians associated with ``orientable string modules". As an application we explicitly compute the Euler--Poincar\'e characteristic of quiver Grassmannians associated with  indecomposable preprojective, preinjective and regular homogeneous representations of an affine quiver of type $\tilde{A}_{p,1}$.  For $p=1$, this approach provides another proof of a result due to P.~Caldero and A.~Zelevinsky in \cite{CZ}.  
\end{abstract}

\maketitle

\subsection*{keywords} cluster algebras, cluster character, quiver Grassmannians, Euler characteristic, string modules.

\section{Introduction and main results}\label{intro}
In this paper we provide a technique to compute the Euler--Poincar\'e characteristic of some complex projective varieties called  \emph{quiver Grassmannians}. In the last few years many authors have shown the importance of such projective varieties and of their Euler--Poincar\'e characteristic in the theory of cluster algebras (see \cite{CC}, \cite{CK1}, \cite{CK2}, \cite{DWZII}), introduced and studied by S.~Fomin and A.~Zelevinsky (\cite{FZI}, \cite{FZII}, \cite{FZIV}). 

Given a quiver $Q$ and a $Q$--representation $M$, the quiver Grassmannian $Gr_\mathbf{e}(M)$ is the set of all sub--representations of $M$ of a fixed dimension vector $\mathbf{e}$ (see section~\ref{Sec:Quivers}).  This is a complex projective variety and our aim is  to compute its Euler--Poincar\'e characterisitc $\chi_\mathbf{e}(M)$. Our main result (theorem~\ref{MainThmIntro}) says that under some technical hypotheses on $M$, there is an algebraic action of the one--dimensional torus $T=\mathbb{C}^\ast$ on $Gr_\mathbf{e}(M)$. It is well--known (see section~\ref{Sec:ProofMainThm}) that if a complex projective variety is endowed with an algebraic action of a complex torus with finitely many fixed points, then its Euler--Poincar\'e characteristic equals the number of fixed points of this action and,  in particular, it is positive. In general it is not true that the Euler--Poincar\'e characteristic of a quiver Grassmannian is positive (see \cite[example~3.6]{DWZII}) but it is proved in \cite{Nakajima} for quiver Grassmannians associated with rigid representations of acyclic quivers,  as conjectured in \cite{FZI}.  The fixed points of the action of $T$ on $Gr_\mathbf{e}(M)$ are the ``coordinate'' subrepresentations of $M$ of dimension vector $\mathbf{e}$ (section~\ref{Sec:MainResult}).  As a combinatorial tool to count them, we consider the coefficient quiver introduced by Ringel (see section~\ref{Sec:CoeffQuiver}) and we notice that its successor closed subquivers are in bijection with coordinate subrepresentations of $M$ (proposition~\ref{Prop:CoeffQuiver}).  

We prove that ``orientable string modules" (see definition~\ref{Def:OrientableString}) satisfy the hypotheses of theorem~\ref{MainThmIntro}. Such a class of $Q$--representations includes (up to ``right--equivalence") all the representations of the affine quiver of type $\tilde{A}_{p,1}$ and most of the representations of the affine quiver of type $\tilde{A}_{p,q}$. 

As an application we explicitly compute  $\chi_\mathbf{e}(M)$  when $M$ is an indecomposable preprojective, preinjective and regular homogeneous representation of the affine quiver of type $\tilde{A}_{p,1}$. We hence find another proof of results of \cite{CZ} for $p=1$, and of \cite{Thesis} and \cite{irelli-2009} for $p=2$. Such computations can be used to have an explicit description of the bases of cluster algebras of type $\tilde{A}_{p,q}$ found in \cite{irelli-2009} and \cite{dupont-2008} and for further studies of such  cluster algebras \cite{CEsp}. In addition it would be interesting to compare our computations with results of \cite{MS} where the authors compute the Laurent expansion of cluster variables of cluster algebras arising from surfaces. In particular this gives a technique to compute the Euler--Poincar\'e characteristic of quiver Grassmannians associated with rigid representations of quivers associated with triangulations of surfaces with marked points. This family includes quivers of type $\tilde{A}_{p,q}$ where our technique applies. In type $A$ one can compare our results with results of  \cite{ARS}. 

To conclude the introduction we remark that having a torus action on a smooth projective variety $X$ gives rise to a cellular decomposition of $X$ (\cite{BiaBi}, \cite{ChrissGinzburg}).  It is known that if $M$ is a rigid $Q$--representation (i.e. without self--extensions) then $Gr_\mathbf{e}(M)$ is smooth \cite{CR}. In particular if $M$ is a rigid $Q$--representation satisfying hypothesis of theorem~\ref{MainThmIntro} then $Gr_\mathbf{e}(M)$ has a cellular decomposition. This approach is used in \cite{CEsp}. 

The paper is organized as follows: in section~\ref{Sec:Quivers} we recall some basic facts about quivers and quiver Grassmannians; in section~\ref{Sec:MainResult} we state our main result; in section~\ref{Sec:CoeffQuiver} we introduce the coefficient quiver of a $Q$--representation and we show how to use it as a combinatorial tool to apply the main result; in section~\ref{Sec:String} we introduce orientable string modules and we prove that they satisfy the hypotheses of our main theorem; in section~\ref{Sec:TypeAp} we give an explicit application for quivers of type $\tilde{A}_{p,1}$. All the remaining sections are devoted to proofs.

\subsection{Quiver Grassmannians}\label{Sec:Quivers}
We recall the definition of quiver Grassmannians. Given a quiver $Q=(Q_0,Q_1)$, i.e. an oriented graph with vertex set $Q_0=\{1,\cdots, n\}$ and arrow set $Q_1$, a \emph{$Q$--representation} $M$ consists of  a collection of complex  vector spaces $\{M(i),\, i\in Q_0\}$ and a collection of linear maps $\{M(a):M(j)\rightarrow M(i)  ~|~ a:j\rightarrow i\in Q_1\}$. 
\begin{example}
The first column of table~\ref{fig:CoveringRingel} shows some examples of quivers $Q$ and the second one shows an example of a $Q$--representation $M$. 
We denote by $k$ the field of complex numbers.  In the last two rows we use the notation $E_{i,j}$ to denote the linear operator on $k^4$ which sends the $j$--th basis vector  to the $i$--th one and fixes all the others. 
\end{example}
\begin{table}[htbp]
$$
\begin{array}{|c|c|c|c|}
\hline
&&&\\
&Q&M&\tilde{Q}(M)\\
\hline
1&\xymatrix{1\ar^a[r]&2}&\xymatrix{k\ar^{\left[\tiny{\begin{array}{c}1\\1\end{array}}\right]}[r]&k^2}\simeq\xymatrix{k\ar^{\left[\tiny{\begin{array}{c}1\\0\end{array}}\right]}[r]&k^2}&
\xymatrix@R=3pt@M=0pt{\bullet\ar^a[r]\ar_a[rd]&\bullet\\&\bullet}\;\;\;\;\;\xymatrix@R=3pt@M=0pt{\bullet\ar^a[r]&\bullet\\&\bullet} \\
\hline
2&\xymatrix{1&2\ar@<1ex>^a[l]\ar_b[l]} & \xymatrix{k&k\ar@<1ex>^{1}[l]\ar_1[l]} & \xymatrix{\bullet&\bullet\ar@<1ex>^a[l]\ar_b[l]} \\
\hline
3&\xymatrix@R=3pt@C=10pt{&2\ar_{a}[dl]&\\
              1& &3\ar[ul]_{b}\ar[ll]^{c}}&\xymatrix@R=3pt@C=10pt{&k^2\ar_{Id_2}[dl]&\\
              k^2& &k\ar[ul]_{\left[\tiny{\begin{array}{c}1\\0\end{array}}\right]}\ar[ll]^{\left[\tiny{\begin{array}{c}0\\1\end{array}}\right]}}&\xymatrix@R=5pt@C=25pt{&*{\bullet}\ar_{a}[ddl]&\\&*{\bullet}\ar_{a}[ddl]& \\*{\bullet}&&\\
              *{\bullet}& &*{\bullet}\ar[uuul]_{b}\ar[ll]^{c}}\\
\hline
4&\xymatrix{1&2\ar@<1ex>^a[l]\ar_b[l]}&\xymatrix{k^2&k^2\ar@<1ex>^{Id}[l]\ar_{J_2(0)}[l]}&\xymatrix@1@R=8pt{\bullet&\bullet\ar_a[l]\\\bullet&\bullet\ar_(.4){b}[ul]\ar^a[l]}
\\\hline
5&\xymatrix{1\ar@(ur,dr)^b\ar@(ul,dl)_a}&\xymatrix{k^4\ar@(ur,dr)^{E_{32}}\ar@(ul,dl)_{E_{21}+E_{43}}}&\xymatrix{\bullet\ar^a[r]&\bullet\ar^b[r]&\bullet\ar^a[r]&\bullet}\\
&&&\\
\hline
6&\xymatrix{1\ar@(ur,dr)^b\ar@(ul,dl)_a}&\xymatrix{k^4\ar@(ur,dr)^{E_{32}}\ar@(ul,dl)_{E_{21}+E_{34}}}&\xymatrix{\bullet\ar^a[r]&\bullet\ar^b[r]&\bullet&\bullet\ar_a[l]}\\
&&&\\
\hline
\end{array}
$$
\caption{Some $Q$--representations and their coefficient--quiver. In the fourth row, we denote by $J_2(0)$ the $2\times 2$ nilpotent Jordan block. In the last two rows $E_{ij}$ denotes the $4\times 4$ elementary matrix with $1$ in the $ij$--component and zero elsewhere.}
\label{fig:CoveringRingel}
\end{table}
A \emph{subrepresentation} $N$ of $M$ consists of a collection of vector subspaces $N(i)$ of $M(i)$, $i\in Q_0$, such that $M(a)N(j)\subset N(i)$ for every arrow $a:j\rightarrow i$ of $Q$. For example the $Q$--representation $M$ shown in the first line of table~\ref{fig:CoveringRingel} does \emph{not} admit  the $Q$--representation $(\xymatrix@1{k\ar[r]&0})$ as its subpresentation (because the  map $M(a)$ has one--dimensional image) but admits $(\xymatrix@1{0\ar[r]&k})$. 

The \emph{dimension vector} of $M$ is the vector $\mathbf{dim}(M):=(dim_{\mathbb{C}} (M(i)): i\in Q_0)$ where $dim_\mathbb{C}(M(i))$ denotes the complex dimension of the vector space $M(i)$.  For example in table~\ref{fig:CoveringRingel} the dimension vector of $M$ is respectively, from above to below, $(1,2)$, $(1,1)$, $(2,2,1)$, $(2,2)$, $(4)$, $(4)$.

The \emph{path algebra} $kQ$ of $Q$ is the complex vector space with as basis the \emph{paths} of $Q$ (i.e. concatenations of arrows) endowed with the multiplication given by the juxtaposition of paths. It is known (see e.g. \cite{AusReiOsm}) that the category of $Q$--representations is equivalent to the category of $kQ$--modules. In particular every $Q$--representation can be seen as a $kQ$--module and viceversa every $kQ$--module has a natural structure of $Q$--representation. 

Finally, the \emph{quiver Grassmannian} $Gr_\mathbf{e}(M)$ of $M$ of dimension $\mathbf{e}=(e_i: i\in Q_0)$ is defined as the set
of all  the subrepresentations of $M$ of dimension vector $\mathbf{e}$, that is, 
$$
Gr_\mathbf{e}(M):=\{N\subset M: \mathbf{dim}(N)=\mathbf{e}\}.
$$
\begin{example}
For the $Q$--representations $M$ shown in lines~1 and 2  of table~\ref{fig:CoveringRingel} the quiver Grassmannian $Gr_{(1,1)}(M)$ is a point.  
If $M$ is the $Q$--representation of line~3, $Gr_{(1,1,1)}(M)$ is the empty set.
Let $M$ be the $Q$--representation shown in line~4. Here $J_2(0)=E_{12}=\left[\tiny{\begin{array}{cc}0&1\\0&0\end{array}}\right]$ is the $2\times 2$ nilpotent Jordan block which sends the second basis vector to the first one.  We consider the set $Gr_{(1,1)}(M)$ of subrepresentations of $M$ of dimension vector $(1,1)$. This consists of lines in $k^2$ spanned by non--zero vectors $v=(\lambda,\mu)^t\in k^2$ such that $v$ and $J_2(0)v$ are linearly dependent. In other words a line spanned by $v$ is in $Gr_{(1,1)}(M)$ if and only if 
$
\rm{det}\left[\begin{array}{cc}\lambda&\mu\\\mu&0\end{array}\right]=-\mu^2=0.
$
Then $Gr_{(1,1)}(M)$ is a  point which is actually not reduced, indeed the tangent space at this point has dimension one (see e.g. \cite{CEsp}). 

If $M$ is the $Q$--representation shown in line~5 we consider  $Gr_{(1)}(M)$ which consists of the lines of $k^4$ invariant under the linear operators $E_{21}+E_{43}$ and $E_{32}$. It is easy to see that this set consists only of the line spanned by the fourth basis vector. Similarly if $M$ is the $Q$--representation shown in the last row of table~\ref{fig:CoveringRingel}, $Gr_{(1)}(M)$ consists only of one point:  the line spanned by the third basis vector. 
\end{example}
We notice that the quiver Grassmannian $Gr_\mathbf{e}(M)$ is closed inside the product $\prod_{i\in Q_0}Gr_{e_i}(M(i))$, where  $Gr_{e_i}(M(i))$ denotes  the usual Grassmannian  of all vector subspaces of $M(i)$ of dimension $e_i$, which is a projective variety. As a consequence, $Gr_\mathbf{e}(M)$ is a complex projective variety. We denote by $\chi_\mathbf{e}(M)$ its Euler--Poincar\'e characteristic. In the examples shown above  $\chi_\mathbf{e}(M)$ is one if $Gr_\mathbf{e}(M)$ is a (double) point and zero if it is the empty set. 

\subsection{The main result}\label{Sec:MainResult} The following theorem is our main result. 
 
\begin{theorem}\label{MainThmIntro} Let $M$ be a $Q$--representation and for every $i\in Q_0$ let $B(i)$ be a  linear basis of $M(i)$ such that  for every arrow $a:j\rightarrow i$ of $Q$ and every element $b\in B(j)$ there exists an element $b'\in B(i)$ and  $c\in k$ (possibly zero) such that
\begin{equation}\label{HypthTreeModIntro}
M(a)b=cb'.
\end{equation}
Suppose that each $v\in B(i)$ and all its multiples $cv$, $c\in k^\ast$, is assigned a degree $d(cv)=d(v)\in\ZZ$ so that:
\begin{enumerate}
\item[(D1)] for all $i\in Q_0$ all vectors from $B(i)$ have different degrees;
\item[(D2)] for every arrow  $a:j\rightarrow i$ of $Q$, whenever $b_1\neq b_2$ are elements of $B(j)$ such that $M(a)b_1$ and $M(a)b_2$ are non-zero we have:
\begin{equation}
d(M(a)b_1)-d(M(a)b_2)=d(b_1)-d(b_2).
\end{equation}
\end{enumerate}
Then 
\begin{equation}\label{eq:Chie(M)}
\chi_\mathbf{e}(M)=|\{N\in Gr_{\mathbf{e}}(M):\, N(i)\text{ is spanned by a part of }B(i)\}|.
\end{equation}
\end{theorem}
The hypothesis \eqref{HypthTreeModIntro} says that every column and every row of the matrix $M(a)$ contains at most one entry different from zero. 

The hypothesis $(D2)$ can be replaced by saying that every arrow $a$ of $Q$ has a degree $d(a)\in\ZZ$ so that $d(b')=d(b)+d(a)$ whenever $M(a)b=cb'$, for some non--zero coefficient $c\in k$. 

The thesis \eqref{eq:Chie(M)} says that we need to count the number of ``coordinate'' subrepresentations i.e. those $N\in Gr_\mathbf{e}(M)$  whose vector space $N(i)$ is a coordinate subspace in the basis $B(i)$ (i.e. is spanned by elements of $B(i)$). 
\begin{example}
Let $Q$ be the quiver with only one vertex and no--arrows. A $Q$--representation is just a vector space $V$ and the quiver Grassmannians are usual Grassmannians of vector subspaces. Let $\{v_1,\cdots, v_n\}$ be a basis of $V$. We assign degree $d(v_i):=i$ and the hypotheses of theorem~\ref{MainThmIntro} are satisfied. Then, by theorem~\ref{MainThmIntro},  $\chi(Gr_k(V))$ is the number of coordinate vector subspaces (i.e. generated by basis vectors) of $V$ of dimension $k$. We hence find the well--known result: $\chi(Gr_k(V))={n\choose k}$.

Let us give other examples with the help of table~\ref{fig:CoveringRingel}. The $Q$--representations shown in line~1 are isomorphic, but the first one does not satisfy the hypothesis~\eqref{HypthTreeModIntro} and we cannot apply theorem~\ref{MainThmIntro}, while the second one does.

The second line shows an interesting example. The $Q$--representation $M$ of this line is a ``deformation" of $M':=\xymatrix@1{k&k\ar@<0.5ex>_{1}[l]\ar^{0}[l]}$ and they have the same quiver Grassmannians (see lemma~\ref{Lemma:RightKronecker}). These two $Q$--representations are indeed \emph{right--equivalent} in the sense of \cite{DWZ}. Theorem~\ref{MainThmIntro} applies to $M'$ and we can hence compute $\chi_\mathbf{e}(M)$. 

In line~3 of table~\ref{fig:CoveringRingel} we choose $d(a)=d(b):=0$ and $d(c):=1$ and hence the choice of a degree for the generator of the one--dimensional vector space at vertex $3$ determines the choice of a degree for the two basis vectors at vertices $2$ and $3$ and  these two degrees are different. We can hence apply theorem~\ref{MainThmIntro}. 

In line~4 we choose $d(a):=0$ and $d(b):=1$.

In line~5 we choose $d(a)=d(b)=1$.

In line~6 we choose $d(a)=1$ and $d(b)=2$. 
\end{example}

\subsection{Coefficient--quiver}\label{Sec:CoeffQuiver}
In order to compute $\chi_\mathbf{e}(M)$ with the help of theorem~\ref{MainThmIntro} one can use a combinatorial tool called the \emph{coeffcient--quiver} $\tilde{Q}(M,B)$ of $M$ in the basis $B$ (introduced by Ringel  in \cite{RingelExceptional}). Let us recall its definition and show its utility. Let $M$ be a $Q$--representation and $B=\cup_{i\in Q_0}B(i)$ a collection of basis $B(i)$ of $M(i)$. The set $B$ is hence a basis of the vector space $\bigoplus_{i\in Q_0}M(i)$ and we refer to it as a \emph{basis} of $M$. The coefficient--quiver  $\tilde{Q}(M,B)$  is a quiver whose vertices are identified with the elements of $B$; the arrows are defined as follows:  for every arrow $a:j\rightarrow i$ of $Q$ and every element $b\in B(j)$ we expand $M(a)b=\sum c_{b'}b'$ in the basis $B(i)$ of $M(i)$ and we put  an arrow (still denoted by $a$) from $b$ to $b'\in B(i)$ in $\tilde{Q}(M,B)$  if the coefficient $c_{b'}$ of $b'$ in this expansion  is non--zero. Table~\ref{fig:CoveringRingel} shows examples of coefficient--quivers (which are denoted simply by $\tilde{Q}(M)$ since they are in the basis  in which $M$ is presented).

We denote by $T\overrightarrow{\subset}\tilde{Q}(M)$ a \emph{successor closed} subquiver $T$ of $\tilde{Q}(M)$, i.e. a subquiver $T$ such that if $j\in T_0$ is one of its vertices and $a:j\rightarrow i$ is an arrow of $Q$ then $a$ is an arrow of $T$.  

It is easy to see that the following proposition is equivalent to theorem~\ref{MainThmIntro}.

\begin{proposition}\label{Prop:CoeffQuiver}
Let $M$ be a $Q$--representation satisfying hypotheses of theorem~\ref{MainThmIntro}. Then 
\begin{equation}\label{eq:ChiMReformulation}
\chi_{\mathbf{e}}(M)=|\{T\overrightarrow{\subset}\tilde{Q}(M)~:~|T_0\cap B(i)|=e_i,\, \forall i\in Q_0\}|
\end{equation} 
where $T_0$ denotes the vertices of $T$. In particular $\chi_{\mathbf{e}}(M)$ is positive.
\end{proposition}
For example let us consider the $Q$--representation $M$ shown in the third line of table~\ref{fig:CoveringRingel}. We have already noticed that $M$ satisfies hypotheses of theorem~\ref{MainThmIntro}. Then we apply proposition~\ref{Prop:CoeffQuiver} and we find $\chi_{(1,0,0)}(M)=2$. Indeed there are two successor closed subquivers of $\tilde{Q}(M)$ with $|T_0\cap B(1)|=2$ and $|T_0\cap B(2)|=|T_0\cap B(3)|=0$ which are the two sinks (this is consistent with the fact that $Gr_{(1,0,0)}(M)=\mathbb{P}^1(k^2)$ is a projective line). Many other examples can be taken from table~\ref{fig:CoveringRingel}.
 
\subsection{String--modules}\label{Sec:String} We now show a class of $Q$--representations which satisfy the hypotheses of theorem~\ref{MainThmIntro}.

A $Q$--representation $M$ is called a \emph{string module} if it admits a basis $B_0$ such that the coefficient--quiver $\tilde{Q}(M,B_0)$ in this basis is a chain (i.e. a  $2$--regular graph not necessarily connected) and if every column and every row of every matrix $M(a)$ in this basis $B_0$ has at most one non--zero entry, i.e.  it satisfies \eqref{HypthTreeModIntro}. We remark that this definition follows \cite{CrawleyTree} but not \cite{RingelExceptional} where \eqref{HypthTreeModIntro} is not required. For a string module $M$ we sometimes avoid mentioning the basis $B_0$ and we denote the corresponding coefficient--quiver simply by $\tilde{Q}(M)$.  The $Q$--representations shown in  table~\ref{fig:CoveringRingel}  are all string modules except the second one. It can be shown that a string module $M$ is indecomposable if and only if $\tilde{Q}(M)$ is connected (\cite{CrawleyTree}, \cite[Sec.~3.5 and 4.1]{GabrielCover}). 

Given an indecomposable string module $M$,  the chain $\tilde{Q}(M)$ has  two extreme vertices (i.e. joined with exactly one vertex). We say that two arrows of $\tilde{Q}(M)$ have the \emph{same orientation} if they both point toward the same extreme vertex and they have \emph{different orientation} otherwise. For example the two arrows labelled by $a$ in  lines 5 and 6 of  table~\ref{fig:CoveringRingel} have the same orientation in the line~5 while they have different orientation in the line~6.  

During private conversations with J.~Schr\"oer we were introduced to the following definition.
\begin{definition}\label{Def:OrientableString}
A string module $M$ is called \emph{orientable} if for every arrow $a$ of $Q$, all the corresponding arrows $a$ of $\tilde{Q}(M)$  have the same orientation.
\end{definition}
For example line~5 of table~\ref{fig:CoveringRingel} shows an orientable string module while the line~6 shows a non--orientable one. 

\begin{proposition}\label{Prop:StringModules}
If $M$ is an orientable string module then \eqref{eq:ChiMReformulation} holds.
\end{proposition} 
In section~\ref{Prop:StringModules} we show that an orientable string module satisfies \eqref{HypthTreeModIntro}, $(D1)$ and $(D2)$ and hence, by  proposition~\ref{Prop:CoeffQuiver}, they satisfy \eqref{eq:ChiMReformulation}.

\subsection{Explicit computations in type $\tilde{A}_{p,1}$}\label{Sec:TypeAp} In this section 
we compute explicitly $\chi_\mathbf{e}(M)$ for some indecomposable representation $M$ of the affine quiver $Q_{p,1}$ of type $\tilde{A}_{p,1}$. Let us recall the definition of $Q_{p,1}$. 

Let  $p\geq1$ be an integer. By definition $Q_{p,1}$ has one sink, one source and $p+1$ arrows which form two paths, one with $p$ arrows and the other with one arrow. We denote the vertices of $Q_{p,1}$ by numbers from $1$ to $p+1$ so that $1$ is the sink, $p+1$ is the source and $k$ is joined to $k+1$ by the arrow $\varepsilon_{k}$, for $k=1,2,\cdots,p$ and $p+1$ is joined to $1$ by the arrow $\varepsilon_0$  as shown below: 
$$
\xymatrix@R=9pt@C=10pt{
                                                 &2\ar_{\varepsilon_1}[d]&  \ar_{\varepsilon_2}[l]\cdots  &p\ar_{\varepsilon_{p-1}}[l]  \\
Q_{p,1}:=            &1        &                & p+1\ar^{\varepsilon_{0}}[ll]\ar_{\varepsilon_{p}}[u]
}
$$
For every $n\geq0$ and $1\leq t\leq p$ we define the $Q_{p,1}$--representations
$$
\begin{array}{cc}
\xymatrix@R=7pt@C=8pt{
&*+[F]{k^{n+1}}\ar[d]&\\
&k^{n+1}\ar[d]&k^{n}\ar_{\varphi_1}[ul]\\
&\vdots\ar[d]&\vdots\ar[u]\\
M^n_p(\overline{[1,t]}):=&k^{n+1}&k^n\ar[u]\ar^{\varphi_2}[l]
}&
\xymatrix@R=7pt{
&*+[F]{k^{n}}\ar[d]&\\
&k^{n}\ar[d]&k^{n+1}\ar_{\varphi_2^t}[ul]\\
&\vdots\ar[d]&\vdots\ar[u]\\
M^n_p(\underline{[1,t]}):=&k^{n}&k^{n+1}\ar[u]\ar^{\varphi_1^t}[l]
}
\end{array}
$$
where the highlighted vector spaces correspond to the vertex $t$. These representations are called respectively pre--projective and pre--injective modules (see e.g.\cite{ASS}).

For every $\lambda\in k$ and $n\geq1$, let $Reg_p^n(\lambda)$ be the $Q_{p,1}$--representation 
$$
\xymatrix@R=7pt@C=8pt{
                                                 &k^{n}\ar_{=}[d]&  \ar_=[l]\cdots  &k^{n}\ar_=[l]  \\
Reg^n_p(\lambda):=            &k^{n}        &                & k^n\ar^{J_n(\lambda)}[ll]\ar_{=}[u]
}
$$
with a Jordan block $J_n(\lambda)$ of eigenvalue $\lambda$ at the arrow $\varepsilon_0$ and the identity map in all the other arrows. This reprsentation is called regular homogeneous. 
It is easy to see that $M_p^n(\overline{[1,t]})$, $M_p^n(\underline{[1,t]})$ and $Reg_p^n(0)$ are orientable string modules (see lemma~\ref{Lemma:Ap1Orientable}) and $\chi_\mathbf{e}(Reg_p^n(\lambda))=\chi_\mathbf{e}(Reg_p^n(0))$ for every  $\lambda\in k$ (section~\ref{Sec:Ap}). 
We can hence apply  theorem~\ref{MainThmIntro} (or proposition~\ref{Prop:StringModules}).

We often use the following notation: 
\begin{equation}\label{eq:Chi[rs]}
\chi_{\mathbf{e}}([r,s]):=\prod_{k=r}^{s-2}{e_k-e_{s}\choose e_{k+1}-e_{s}}=\prod_{k=r+1}^{s-1}{e_r-e_{k+1}\choose e_{k}-e_{k+1}}
\end{equation}
with the convention that this product equals one whenever $r> s-2$. We interpret $\chi_{\mathbf{e}}([r,s])$ as the Euler characteristic of the flag variety 
$$
\{k^{e_r}\supseteq M_{r+1}\supseteq\cdots\supseteq M_{s-1}\supseteq k^{e_s}|\;dim (M_k)=e_k\}.
$$
\begin{proposition}\label{Prop:Ap1}
For every $n\geq1$, $1\leq t\leq p$ and $\lambda\in k$ we have

\begin{equation}\label{Eq:M1t+}
\chi_{(e_1,\cdots,e_{p+1})}(M^n_p(\overline{[1,t]}))=
\end{equation}
$$
{e_1-1\choose e_{p+1}}{n\!+\!1\!-\!e_{t}\choose e_1-e_{t}}{n\!+\!1\!-\!e_{t+1}\choose e_{t}-e_{t+1}}{n\!-\!e_{p+1}\choose e_{t+1}-e_{p+1}}
\chi_{\mathbf{e}}([1,t])\chi_{\mathbf{e}}([t+1,p+1])
$$
\begin{equation}\label{Eq:M1t-}
\chi_{(e_1,\cdots,e_{p+1})}(M^n_p(\underline{[1,t]}))=
\end{equation}
$$
{n-e_{p+1}\choose e_1-e_{p+1}}{e_{t+1}\choose e_{p+1}}{e_{t}+1\choose e_{t+1}}{e_1\choose e_{t}}
\chi_{\mathbf{e}}([1,t])\chi_{\mathbf{e}}([t+1,p+1])
$$
\begin{equation}\label{eq:ChiRegular}
\chi_\mathbf{e}(Reg^n_{p}(\lambda))={e_1\choose e_{p+1}}{n-e_{p+1}\choose e_1-e_{p+1}}\chi_{\mathbf{e}}([1,p+1])
\end{equation}
We always use the convention that the binomial coefficient ${p\choose q}$ equals 0 if $q<0$, $p<0$, $q>p$ and it equals 1 if $q=0$ and $p\geq q$.
\end{proposition}
 
\section{Proof of theorem~\ref{MainThmIntro}}\label{Sec:ProofMainThm}

The proof is based on the following well--known fact:
given a complex projective variety $X$ and an algebraic action $\varphi:T\times X\rightarrow X$, $(\lambda,x)\mapsto \lambda.x$ of the one--dimensional torus $T=\mathbb{C}^\ast$ with finitely many fixed points, then the number of fixed points equals the Euler--Poincar\'e characteristic $\chi(X)$ of $X$. To see this we consider the decomposition $X=X^T\coprod Y$ of $X$ into the disjoint union of the set $X^T$ of  fixed points of $\varphi$ and of their complement $Y:=X\setminus X^T$. Such sets are locally closed and hence $\chi(X)=\chi(X^T)+\chi(Y)$. The restriction of $\varphi$ to $Y$ defines a surjective morphism $\varphi: T\times Y\rightarrow Y$ whose fibers are all isomorphic to $\mathbb{C}^\ast$. It follows that $\chi(Y)=\chi(\mathbb{C}^\ast)=0$ and hence $\chi(X)=\chi(X^T)$ which equals the number of fixed points of $\varphi$.

We hence find a torus action on  our quiver Grassmannians. 

Let M be a representation satisfying hypotheses (D1) and (D2) of the theorem.  The torus $k^{*}$ acts on $M$  as follows: 
\begin{eqnarray}\label{eq:TorusAction}
\lambda.b:=\lambda^{d(b)}b,&&\lambda\in \emph{k}^{*}
\end{eqnarray}
for every element $b\in B$ of the basis $B$ extended by linearity to all the elements of $M$.  This action extends to quiver Grassmannians:

\begin{lemma}\label{lemma:ActionOfTorus} 
Let $U\in Gr_{\mathbf{e}}(M)$ be a subrepresention of $M$ of dimension vector $\mathbf{e}$. Then, given $\lambda\in \emph{k}^{*}$, the set $\lambda.U:=\{\lambda.u|\,u\in U\}$ is a subrepresentation of $M$ of the same dimension vector $\mathbf{e}$ of $U$.  
\end{lemma}
\begin{proof}
Given an arrow $a:j\rightarrow i$ of $Q$ we define the number  $d(a):=d(M(a)b)-d(b)$ for an element $b\in B(j)$ such that $M(a)b$ is non--zero. This definition is independent of the choice of $b$ in view of (D2). Then it is easy to verify that for every $v\in M(j)$
$$
\lambda. (M(a)v)=\lambda^{d(a)}M(a)(\lambda. v)
$$
which concludes the proof.
\end{proof}

Given a subrepresentation $U\in Gr_\mathbf{e}(M)$,  the element $\lambda\in k^*$ acts on each vector subspace $U(i)$ as a diagonal operator with different eigenvalues, in view of property (D1). Then the fixed subrepresentations $U=\lambda.U\in Gr_\mathbf{e}(M)$ are precisely the coordinate subspaces of $M$ in the basis $B$ of dimension $e:=\sum_i e_i$ which concludes the proof of theorem~\ref{MainThmIntro}.

\section{Proof of proposition~\ref{Prop:StringModules}}
We prove that an orientable string module $M$ satisfies the hypotheses of theorem~\ref{MainThmIntro}.
By definition there exists a basis $B_0$ of $M$ so that \eqref{HypthTreeModIntro} is satisfied and the coefficient--quiver $\tilde{Q}(M,B_0)$ in $B_0$ is a chain. We have to assign a degree $d(b)\in\ZZ$ to the elements of $B_0$ (which are also the vertices of $\tilde{Q}(M,B_0)$) so that (D1) and (D2) are satisfied. 

Since $S:=\tilde{Q}(M,B_0)$ is a chain we number  the vertices of $S$ as $s_1,s_2,\cdots$ in such a way that for every $i=1,\cdots,m$ there is a unique edge $\varepsilon_{i}$ between $s_i$ and $s_{i+1}$.   We assign the degree $d(s_i):=i$ for $i=1,2, \cdots$. Then (D1) is clearly satisfied (all the elements of $B_0$ have different degrees and hence all the elements of $B_0(i)$ have different degrees). Since $M$ is orientable it is also easy to prove that (D2) is satisfied. Indeed, by definition, for every arrow $a$ of $Q$ all the corresponding arrows $a$ of $S$ have all the same orientation, either all of them are oriented from $s_i$ to $s_{i+1}$ or from $s_{i+1}$ to $s_i$.

\section{Proof of proposition~\ref{Prop:Ap1}}

For the convenience of the reader we prove proposition~\ref{Prop:Ap1} first in the case $p=1$ (the Kronecker quiver) and hence for   $p\geq1$. 

All the proofs are based on the following lemma.
\begin{lemma}\label{Lemma:Ap1Orientable}
$M_p^n(\overline{[1,t]})$, $M_p^n(\underline{[1,t]})$ and $Reg_p^n(0)$ are orientable string modules (in the sense of definition~\ref{Def:OrientableString}). In particular \eqref{eq:ChiMReformulation} holds.
\end{lemma}
\begin{proof}
All the linear maps defining such $Q_{p,1}$--representations satisfy \eqref{HypthTreeModIntro}. It remains to show that their coefficient--quiver is a chain. 
 
Let $S_{\varepsilon_0}$ be the subquiver of  $Q_{p,1}$ obtained by removing the arrow $\varepsilon_0$. We join together $n$ copies of $S_{\varepsilon_0}$ by using the arrow $\varepsilon_0$ and we get a string  that we denote by $S^n_0$. The coefficient--quiver of $Reg_p^n(0)$ is $S^n_0$ which is a chain. 

Let $1\leq t\leq p$ be a vertex of $Q_{p,1}$. We consider the full subquiver $S(\overline{[1,t]})$  of $Q_{p,1}$ with vertex set all the vertices $1,2,\cdots,t$.  We join the string $S^n_0$ with the string $S([1,t])$ by using the arrow $\varepsilon_{0}$ and we get a new string that we call $S^n(\overline{[1,t]})$. Such a string is the coefficient--quiver of $M_{p}^n(\overline{[1,t]})$

In order to get the coefficient--quiver of $M_{p}^n(\underline{[1,t]})$ we proceed similarly: we consider the full subquiver $S(\underline{[1,t]})$ with vertices $t+1,t+2,\cdots, p, p+1$. We join $S(\underline{[1,t]})$ with $S^n$ by using the arrow $\varepsilon_0$ and we get a quiver $S^n(\underline{[1,t]})$. Such a quiver is the coefficient--quiver of $M_{p}^n(\underline{[1,t]})$. Figure~\ref{fig:CoveringM43[13]} shows the case $p=4$, $t=n=3$.

\begin{figure}[htbp]
\begin{center}
$$
\xymatrix@R=6pt@C=0pt{
&             &                       &                      &              &*{5}\ar_{\varepsilon_4}[dl]\ar[dddd]^{\varepsilon_0}&                     &                    &                        &5\ar_{\varepsilon_4}[dl]\ar[dddd]^{\varepsilon_0}&             &                    &                     &*{5}\ar_{\varepsilon_4}[dl]&                      &\\
   &          &                       &                      &4\ar^{\varepsilon_3}[dl]   &                                        &                      &                     &4\ar^{\varepsilon_3}[dl]&                                   &              &                    & 4\ar^{\varepsilon_3}[dl]&                                      &                     &\\
      &       &                       &3\ar^{\varepsilon_{2}}[dl] &              &                                        &                      &3\ar^{\varepsilon_{2}}[dl]&                     &                                   &             &3\ar^{\varepsilon_{2}}[dl]&                      &                                      &                     &\\
         &    &2 \ar^{\varepsilon_1}[dl]           &                      &              &                                        &2\ar^{\varepsilon_1}[dl] &                     &                     &                                    &2\ar^{\varepsilon_1}[dl] &                     &                     &                                      &&                           \\
S^3_0=&1&                       &                      &              & 1                           &                      &                     &                    &           1        &              &                     &                     &             &                      &
}
$$
$$
\xymatrix@R=6pt@C=0pt{
&             &                       &                      &              &*{5}\ar_{\varepsilon_4}[dl]\ar[dddd]^{\varepsilon_0}&                     &                    &                        &5\ar_{\varepsilon_4}[dl]\ar[dddd]^{\varepsilon_0}&             &                    &                     &*{5}\ar_{\varepsilon_4}[dl] \ar[dddd]^{\varepsilon_0} &                      &\\
   &          &                       &                      &4\ar^{\varepsilon_3}[dl]   &                                        &                      &                     &4\ar^{\varepsilon_3}[dl]&                                   &              &                    & 4\ar^{\varepsilon_3}[dl]&                                      &                     &\\
      &       &                       &3\ar^{\varepsilon_{2}}[dl] &              &                                        &                      &3\ar^{\varepsilon_{2}}[dl]&                     &                                   &             &3\ar^{\varepsilon_{2}}[dl]&                      &                                      &                     &3\ar^{\varepsilon_{2}}[dl]\\
         &    &2 \ar^{\varepsilon_1}[dl]           &                      &              &                                        &2\ar^{\varepsilon_1}[dl] &                     &                     &                                    &2\ar^{\varepsilon_1}[dl] &                     &                     &                                      &2\ar^{\varepsilon_1}[dl]&                           \\
S^3(\overline{[1,3]})=&1&                       &                      &              & 1                           &                      &                     &                    &           1        &              &                     &                     &1                 &                      &
}
$$
$$
\xymatrix@R=6pt@C=0pt{
&             &                       &                      &              &*{5}\ar_{\varepsilon_4}[dl]\ar[dddd]^{\varepsilon_0}&                     &                    &                        &5\ar_{\varepsilon_4}[dl]\ar[dddd]^{\varepsilon_0}&             &                    &                     &*{5}\ar_{\varepsilon_4}[dl] \ar[dddd]^{\varepsilon_0} &                      &&&*{5}\ar_{\varepsilon_4}[dl]\\
   &          &                       &                      &4   &                                        &                      &                     &4\ar^{\varepsilon_3}[dl]&                                   &              &                    & 4\ar^{\varepsilon_3}[dl]&                                      &                     & &4\ar^{\varepsilon_3}[dl]\\
      &       &                       & &              &                                        &                      &3\ar^{\varepsilon_{2}}[dl]&                     &                                   &             &3\ar^{\varepsilon_{2}}[dl]&                      &                                      &                     &3\ar^{\varepsilon_{2}}[dl]\\
         &    &            &                      &              &                                        &2\ar^{\varepsilon_1}[dl] &                     &                     &                                    &2\ar^{\varepsilon_1}[dl] &                     &                     &                                      &2\ar^{\varepsilon_1}[dl]&                           \\
S^3(\underline{[1,3]})=& &                     &                      &              & 1                           &                      &                     &                    &           1        &              &                     &                     &1                 &                      &
}
$$
\caption{The coefficient--quiver of $Reg_4^3(0)$, $M_4^3(\overline{[1,3]})$ and $M_4^3(\underline{[1,3]})$ respectively}
\label{fig:CoveringM43[13]}
\end{center}
\end{figure}
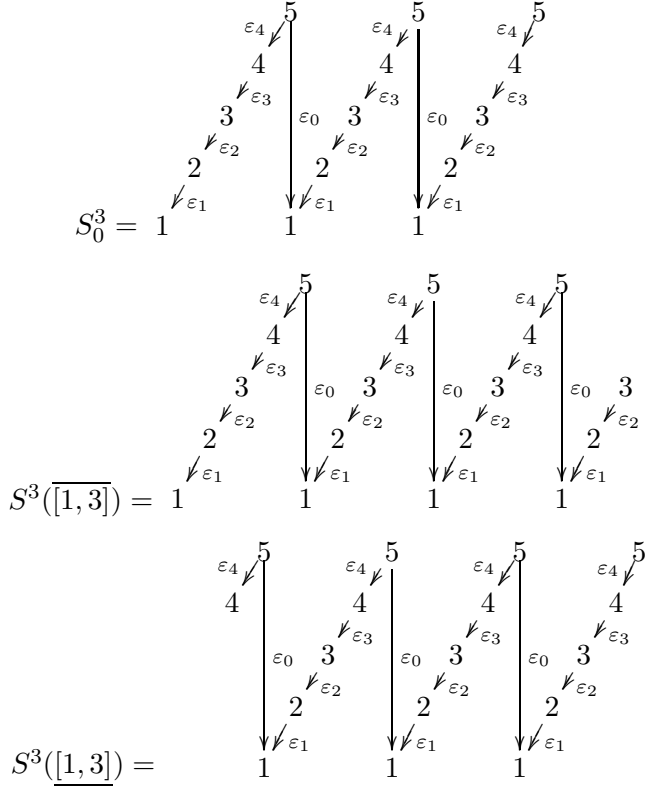
\end{proof}

\subsection{Type $\tilde{A}_{1,1}$: the Kronecker quiver}\label{section:Kronecker}
In this section we consider the Kronecker quiver
$Q_{1,1}:=\xymatrix@1{1&2\ar@<1ex>[l]^{\varepsilon_0}\ar[l]_{\varepsilon_1}}$ and its representations over the field $k$ of complex numbers.  Let $\varphi_1,\varphi_2:k^n\rightarrow k^{n+1}$ be respectively the immersion in the vector subspace spanned respectively by the first and by the last $n$ basis vectors. For every $n\geq0$ and $\lambda\in k$ we consider the representations
$$
\begin{array}{cc}
M_1^n(\overline{[1,1]}):=\xymatrix@1{k^{n+1}&k^{n}\ar@<1ex>[l]^{\varphi_{2}}\ar[l]_{\varphi_{1}}};&
M_1^{n}{\underline{[1,1]}}:=\xymatrix@1{k^{n}&k^{n+1}\ar@<1ex>[l]^{\varphi_{2}^{t}}\ar[l]_{\varphi_{1}^{t}}}\\Reg_1^n(\lambda):=\xymatrix{k^{n}&k^{n}\ar@<1ex>[l]^{J_{n}(\lambda)}\ar[l]_{=}.}&\end{array}$$
The next result is contained in \cite{CZ}. We give a slightly different  proof by using theorem~\ref{MainThmIntro}. 
\begin{proposition}\label{prop:EulerKronecker}\cite[Propositions~4.3 and 5.3]{CZ} 
For every dimension vector $\mathbf{e}=(e_1,e_2)$ and $n\geq0$ we have:
\begin{eqnarray}\label{eq:EulerKroneckerMn+3}
\chi_{(e_{1},e_{2})}(M_1^n(\overline{[1,1]}))&=&{n+1-e_{2}\choose n+1-e_{1}}{e_{1}-1\choose e_{2}}+\delta_{e_{1},0}\delta_{e_{2},0}\\\label{eq:EulerKroneckerM-n}
\chi_{(e_{1},e_{2})}(M_1^n(\underline{[1,1]}))&=&{e_1+1\choose e_2}{n-e_{2}\choose n-e_{1}}+\delta_{e_1,n}\delta_{e_2,n+1}
\end{eqnarray}                                                    
where $\delta_{a,b}$ denotes the Kronecker delta. For every $\lambda\in k$:                                                     
\begin{equation}\label{eq:EulerKroneckerRegular}
\chi_{(e_{1},e_{2})}(Reg_1^n(\lambda))={n-e_{2}\choose n-e_{1}}{e_{1}\choose e_{2}}
\end{equation}
\end{proposition}

\begin{proof}We notice that \eqref{eq:EulerKroneckerM-n} follows from \eqref{eq:EulerKroneckerMn+3}. Indeed $M_1^{n}{\underline{[1,1]}} \simeq DM_1^n(\overline{[1,1]})$ where $D=Hom_k(\cdot,k)$ is the duality functor
and the isomorphism follows by exchanging the two vertices. Then we have (see also \cite[Sec.~1.2]{CR}):
$$
\chi_{(e_{1},e_{2})} (M_1^{n}(\underline{[1,1]}) ) = \chi_{(n+1-e_{2} ,n-e_{1})} (M_1^n(\overline{[1,1]})).
$$
We hence prove \eqref{eq:EulerKroneckerMn+3}. By lemma~\ref{Lemma:Ap1Orientable}, the representation $M_{1}^n(\overline{[1,1]})$ is an orientable string module and we can apply theorem~\ref{MainThmIntro}. In order to compute $\chi_{(e_1,e_2)}(M_{1}^n(\overline{[1,1]}))$, we have hence to count couples $\{T_1,T_2\}$ of subsets $T_1\subset [1,n+1]$, $T_2\subset [1,n]$ such that $|T_i|=e_i$ ($i=1,2$) and $\varphi_1(T_2)\subset T_1$, $\varphi_2(T_2)\subset T_1$ where $\varphi_1,\varphi_2:[1,n]\rightarrow [1,n+1]$ are the two maps defined by $\varphi_1(k)=k$ and $\varphi_2(k)=k+1$ for $k=1,2,\cdots,n$ (here and in the sequel we use the notation $[1,m]:=\{1,2,\cdots,m\}$). 
We need the following lemma. 
\begin{lemma}\cite[proof~of~ proposition~4.3]{CZ}\label{lemma:Combinatorics}
Let $n$ and $r$ be positive integers such that $1\leq r\leq n$. For an $r$--element subset $J$ of $[1,n]$ we denote by $c(J)$ the number of connected components of $J$ (i.e. the number of  maximal connected intervals in $J$). The number of $r$--element subsets $J$ of $[1,n]$ such that $c(J)=c$ is ${r-1\choose c-1}{n+1-r\choose c}$.
\end{lemma}
\begin{proof} A proof of lemma \ref{lemma:Combinatorics} can be found in \cite[proof~of~ proposition~4.3]{CZ}\end{proof} We hence continue the proof of \eqref{eq:EulerKroneckerMn+3}. The choice of an element $k\in [1,n]$ determines the choice of the two different elements $\varphi_1(k)$ and $\varphi_2(k)$ of $[1,n+1]$; in general the choice of a subset $T_2$ of $[1,n]$  of cardinality $e_{2}$ with $c$ connected components determines the choice of 
$c+e_{2}$ elements of $[1,n+1]$. Given such a set $T_2$, there are hence ${n+1-(c+e_{2})\choose e_{1}-(c+e_{2})}$ choices for the sets $T_1$ such that $\{T_1,T_2\}$ is a desired couple. If $e_1=e_2=0$ then $\chi_{(0,0)}(M_1^n(\overline{[1,1]}))=1$. We assume $e_1\geq e_2\geq1$.
By lemma~\ref{lemma:Combinatorics} the number of $e_{2}$--element subsets $T_2$ of $[1,n]$ with $c(T_2)=c$ equals
$
{e_{2}-1\choose c-1}{n+1-e_{2}\choose c}$.
The number of desired couples $\{T_1,T_2\}$ is hence
\begin{eqnarray}
\chi_{(e_{1},e_{2})}(M_1^n(\overline{[1,1]}))&=&\sum_{c=1}^{e_{1}-e_{2}}{n+1-(c+e_{2})\choose e_{1}-(c+e_{2})}{e_{2}-1\choose c-1}{n+1-e_{2}\choose c}\nonumber\\ &=&
\sum_{c=1}^{e_{1}-e_{2}}{e_{1}-e_{2}\choose c}{e_{2}-1\choose c-1}{n+1-e_{2}\choose e_{1}-e_{2}}\nonumber\\
&=&{n+1-e_{2}\choose e_{1}-e_{2}}\sum_{c=1}^{e_{1}-e_{2}}{e_{1}-e_{2}\choose c}{e_{2}-1\choose e_{2}-c}\nonumber\\
&=&{n+1-e_{2}\choose e_{1}-e_{2}}{e_{1}-1\choose e_{2}}\nonumber
\end{eqnarray}
In the second equality we have used the identity:
$
{n+1-r-q\choose p-q}{n+1-r\choose q}={p\choose q}{n+1-r\choose p}
$
with $q=c$, $p=e_{1}-e_{2}$ and $r=e_{2}$; in the last equality we have used the Vandermonde's identity:
$
\sum_{k}{a\choose k}{b\choose c-k}={a+b\choose c}.
$

We now prove \eqref{eq:EulerKroneckerRegular}. We first assume that $\lambda=0$. The representation $Reg^n_{1,1}(0)$ is an orientable string module and we apply theorem~\ref{MainThmIntro}. We prove \eqref{eq:EulerKroneckerRegular} by induction on $n\geq1$. For $n=0$ it is clear. Let hence $n\geq 1$.  We have hence to count the number of couples $\{T_1,T_2\}$ of subsets $T_2\subset T_1\subset[1,n]$ such that $|T_i|=e_i$ and $J_n(0)T_2\subset T_1\cup{0}$ where $J_n(0):[1,n]\rightarrow [1,n]\cup\{0\}$ maps $k$ to $k-1$ for $k=1,2,\cdots,n$. Alternatively, by proposition~\ref{Prop:CoeffQuiver}, we can consider the coefficient quiver $\tilde{Q}(Reg_1^n)$ of $Reg_1^n$  (shown in figure~\ref{fig:CoveringKronecker}) and count its successor closed subquivers with $e_1$ sources and $e_2$ sinks.
\begin{figure}[htbp]
\begin{center}
\xymatrix@M=0pt@R=5pt@C=20pt{
&&&  \bullet           &                           &   &                  & \\
&&& \bullet         &\bullet\ar[ul]\ar[l]                & &    \bullet           & *+[F]{\bullet}\ar[l]            \\  
&&&\bullet         &\bullet\ar[ul]\ar[l]              &   &  \bullet             &\bullet\ar[ul]\ar[l]\\      
&&&\vdots  &    \vdots                    & & \vdots      &    \vdots      \\
&&& \bullet           &                               & & \bullet              &                \\
&&\tilde{Q}(M_1^n(\overline{[1,1]})):&  \bullet       & \bullet;\ar[ul]\ar[l] &\tilde{Q}(Reg_1^n):&  \bullet       & \bullet\ar[ul]\ar[l]\\
&&&       &                           & &                  &            \\
&&&     &                  & &                            &            \\
&&M_1^n(\overline{[1,1]}):&   k^{n+1}&k^{n}\ar@<1ex>[l]^{\varphi_{2}}\ar[l]_{\varphi_{1}};&Reg_1^n:&   k^{n}&k^{n}\ar@<1ex>[l]^{=}\ar[l]_{J_{n}(0)}
}
\caption{Coefficient--quiver of $Q_{1,1}$--representations.}
\label{fig:CoveringKronecker}
\end{center}
\end{figure}
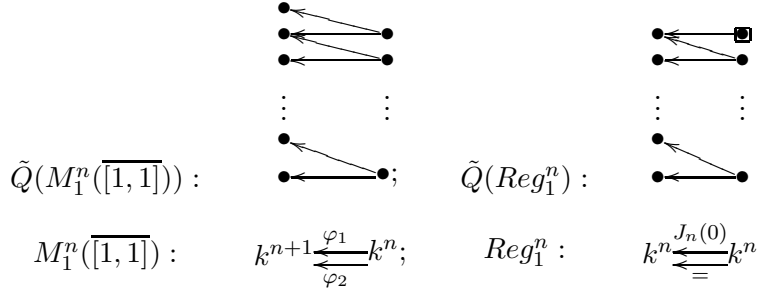
Such a subquiver either contains the unique vertex of $\tilde{Q}(Reg_1^n)$  which is the source of a unique arrow (highlighted in figure~\ref{fig:CoveringKronecker}) or it does not. Alternatively either $T_2$ contains $1=Ker(J_n(0))$ or it does not. We hence have
\begin{eqnarray}\label{Regn1=Mn1}
\chi_{(e_1,e_2)}(Reg_1^n(0))&=&\chi_{(e_1-1,e_2-1)}(Reg^{n-1}_1(0))+\chi_{(e_1,e_2)}(M^{n-1}_1(\overline{[1,1]}))\nonumber\\
&=&{n-e_{2}\choose n-e_{1}}{e_{1}-1\choose e_{2}-1}\!\!+\!\!{n-e_{2}\choose n-e_{1}}{e_{1}-1\choose e_{2}}\!\!+\!\!\delta_{e_{1},0}\delta_{e_{2},0}\nonumber\\
&=&{n-e_{2}\choose n-e_{1}}{e_{1}\choose e_{2}}.\nonumber
\end{eqnarray}
and we are done (we use the obvious fact that ${a-1\choose b-1}+{a-1\choose b}={a\choose b}-\delta_{a,0}\delta_{b,0}$).

It remains to be considered the case where $\lambda\neq0$ which is solved in the following lemma.
\begin{lemma}\label{Lemma:RightKronecker} 
For every $\lambda\in \mathbb{C}$ and $n\geq1$ we have
$$
\chi_\mathbf{e}(Reg_1^n(\lambda))=\chi_\mathbf{e}(Reg_1^n(0))
$$
\end{lemma}
\begin{proof}
As vector spaces, $Reg^n_1(0)$ and $Reg^n_1(\lambda)$ are isomorphic to $k^{2n}$. The path algebra $kQ_{1,1}$ acts on these  isomorphic vector spaces by two actions that we denote respectively by $\ast$ and $\circ$. 
We consider the automorphism $\psi$ of the path algebra $kQ_{1,1}$ which sends $\varepsilon_0$ to $\varepsilon_0+\lambda\varepsilon_1$.  For every $\sigma$ in $kQ_{1,1}$ and every $m$ in $Reg^n_{1,1}(0)$, $\psi(\sigma)\ast m=\sigma\circ m$. Roughly speaking  what the automorphism $\psi$ does is the following: the arrow $\varepsilon_0$ acts as $J_n(0)$ on $Reg_{1,1}^n(0)$, while the arrow $\varepsilon_1$ acts as the identity. Then $\psi(\varepsilon_0)$ acts as $J_n(0)+\lambda Id=J_n(\lambda)$. With this action $Reg_{1,1}^n(0)$ is isomorphic to $Reg_{1,1}^n(\lambda)$ (as $kQ_{11}$--module). In particular the two representations have the same quiver Grassmannians. This proves that they are right--equivalent in the sense of \cite{DWZ}. 
\end{proof}
This concludes the proof of proposition~\ref{prop:EulerKronecker}.
\end{proof}

\subsection{Type $\tilde{A}_{p,1}$}\label{Sec:Ap}
We prove proposition~\ref{Prop:Ap1} for every $p\geq 2$.
The duality functor $D$ sends a representation of $Q_{p,1}$ to a representation of the opposite quiver $Q_{p,1}^{op}$. The symmetries of such quiver induce an isomorphism $M^n_p(\underline{[1,t]})\simeq DM^n_p(\overline{[1,p+1-t]})$ and,  for every dimension vector $\mathbf{e}=(e_1,\cdots,e_{p+1})$, we have:
$$
\chi_{\mathbf{e}}(M^n_p(\underline{[1,t]}))=\chi_{(d_{p+1}-e_{p+1},	\cdots,d_1-e_1)}(M^n_p(\overline{[1,p+1-t]}))
$$
where $\mathbf{d}=(d_1,\cdots,d_{p+1})$ is the dimension vector of $M^n_p(\underline{[1,t]})$. Then \eqref{Eq:M1t-} follows from \eqref{Eq:M1t+}.

We prove \eqref{Eq:M1t+}. By lemma~\ref{Lemma:Ap1Orientable}, the representation $M_p^n(\overline{[1,t]})$ satisfies the hypotheses of theorem~\ref{MainThmIntro}.
In order to compute $\chi_{\mathbf{e}}(M^n_p(\overline{[1,t]}))$ we hence have to count sets $\{T_1,\cdots, T_{p+1}\}$ of subsets $T_1,\cdots, T_t\subset[1,n+1]$, $T_{t+1},\cdots, T_{p+1}\subset[1,n]$ such that: 
$|T_i|=e_i$ and $\varphi_1(T_{t+1})\subset T_t$, $\varphi_2(T_{p+1})\subset T_1$ and $T_k\subset T_{k-1}$ ($k\neq t+1$, $k\neq p+1$) where $\varphi_1, \varphi_2:[1,n]\rightarrow [1,n+1]$ are defined by $\varphi_1(k):=k$ and $\varphi_2(k):=k+1$ for every $k=1,\cdots n$.

For a choice of the quadruple $\{T_1,T_t,T_{t+1},T_{p+1}\}$ (this set could collapse to a  quadruple in which two elements coincide but it does not make any difference in the sequel and we still refer to it as a quadruple) there are $\chi_{\mathbf{e}}([1,t])$ choices for $\{T_2,\cdots,T_{t-1}\}$ and  $\chi_\mathbf{e}([t+1,p+1])$ choices for $\{T_{t+2},\cdots,T_p\}$ such that $\{T_1,\cdots, T_{p+1}\}$ is a desired tuple. 

We hence prove that the number of quadruples $\{T_1,T_t,T_{t+1},T_{p+1}\}$ equals:
\begin{equation}\label{M4}
{e_1-1\choose e_{p+1}}{n\!+\!1\!-\!e_{t}\choose e_1-e_{t}}{n\!+\!1\!-\!e_{t+1}\choose e_{t}-e_{t+1}}{n\!-\!e_{p+1}\choose e_{t+1}-e_{p+1}}
\end{equation}
from which \eqref{Eq:M1t+} follows. 
We hence have to count the number of quadruples $\{T_1,T_{t},T_{t+1},T_{p+1}\}$ of subsets $T_{t}\subset T_1\subset [1,n+1]$, $T_{p+1}\subset T_{t+1}\subset [1,n]$ such that $|T_i|=e_i$, $\varphi_1(T_{t+1})\subset T_{t} $ and $\varphi_2(T_{p+1})\subset T_1$.  

We need the following lemma.
\begin{lemma}\label{lemma:Combinatorics2}
Let $n$ and $e$ be positive integers such that $1\leq e\leq n$. As before, we denote by $c(J)$ the number of connected components of  an $e$--element subset $J$ of $[1,n]$. For every integer $c$, we have 
\begin{enumerate}
\item the number of $e$--element subsets $J$ of $[1,n]$ such that $c(J)=c$ and $J$ contains $n$, is ${e-1\choose c-1}{n-e\choose c-1}$;
\item the number of $e$--element subsets $J$ of $[1,n]$ such that $c(J)=c$ and $J$ does not contain $n$ is ${e-1\choose c-1}{n-e\choose c}$;
\item for every $0\leq r\leq q\leq p$, ${p\choose q}{q\choose r}={p\choose r}{p-r\choose q-r}$. 
\end{enumerate}
\end{lemma}
\begin{proof}
The proof of lemma~\ref{lemma:Combinatorics2} follows from lemma~\ref{lemma:Combinatorics} by an easy induction.
\end{proof}
$$
\begin{array}{ccc}
\xymatrix@R=7pt@C=8pt{
&k^{n+1}\ar[d]&\\
&k^{n+1}\ar[d]&k^{n}\ar_{\varphi_1}[ul]\\
&\vdots\ar[d]&\vdots\ar[u]\\
M^n_p(\overline{[1,t]}):=&k^{n+1}&k^n\ar[u]\ar^{\varphi_2}[l]
}&
\xymatrix@R=7pt@!C{
[1,n+1]\supset\!\!\!\!\!\!\!\!\!\!\!\!\!\!\!\!\!&T_{t}\ar@{_{(}->}[dddd]&\\
                           &                     &T_{t+1}\ar@{_{(}->}_{\varphi_1}[ul]&\!\!\!\!\!\!\!\!\!\!\!\!\!\!\!\!\!\!\!\!\!\subset [1,n]\\
                           &&\\
                           &&\\
[1,n+1]\supset\!\!\!\!\!\!\!\!\!\!\!\!\!\!\!\!\!&T_1&T_{p+1}\ar@{_{(}->}[uuu]\ar@{_{(}->}^{\varphi_2}[l]&\!\!\!\!\!\!\!\!\!\!\!\!\!\!\!\!\!\!\!\!\subset [1,n]
}
\end{array}
$$

Let $T_{p+1}$ be an $e_{p+1}$--element subset of $[1,n]$ and let us count the number of desired quadruples $\{T_1,T_{t},T_{t+1},T_{p+1}\}$ containing $T_{p+1}$. We notice that $T_1$ contains both $\varphi_1(T_{p+1})$ and $\varphi_2(T_{p+1})$. In particular, if $c$ denotes  the number of connected components of $T_{p+1}$, then $T_1$ must contain $c+e_{p+1}$ elements of $[1,n+1]$.  We distinguish the two cases: either $T_{p+1}$ contains $n$ or it does not. 
\begin{enumerate}
\item If $T_{p+1}$ contains $n$ (by lemma~\ref{lemma:Combinatorics2} there are ${e_{p+1}-1\choose c-1}{n-e_{p+1}\choose c-1}$ choices for such subsets) 
then  every possible  $T_1$  contains the element $\varphi_2(n)=(n+1)$.  Then there are ${n+1-c-e_{p+1}\choose e_1-c-e_{p+1}}$ choices for $T_1$. Now either $T_{t}$ contains $(n+1)$ or it does not. If it contains $(n+1)$ (there are ${e_1-1-e_{p+1}\choose e_{t}-1-e_{p+1}}$ choices for such sets) then there are ${e_{t}-1-e_{p+1}\choose e_{t+1}-e_{p+1}}$ choices for $T_{t+1}$; if $T_{t}$ does not contain $(n+1)$ (there are ${e_1-1-e_{p+1}\choose e_{t}-e_{p+1}}$ choices for such sets) then there are ${e_{t}-e_{p+1}\choose e_{t+1}-e_{p+1}}$ choices for $T_{t+1}$. 

The number of quadruples $\{T_1,T_{t},T_{t+1},T_{p+1}\}$ such that $T_{p+1}$ contains $n$ is hence given by:
\begin{eqnarray}\nonumber
\sum_c{e_{p+1}-1\choose c-1}{n-e_{p+1}\choose c-1}{n-e_{p+1}-(c-1)\choose e_1-e_{p+1}-c}\cdot&&\\\nonumber	
[{e_1-1-e_{p+1}\choose e_{t}-1-e_{p+1}}{e_{t}-1-e_{p+1}\choose e_{t+1}-e_{p+1}}+{e_1-1-e_{p+1}\choose e_{t}-e_{p+1}}{e_{t}-e_{p+1}\choose e_{t+1}-e_{p+1}}]&&=\\\nonumber	
\sum_c{e_{p+1}-1\choose c-1}{n-e_{p+1}\choose c-1}{n-e_{p+1}-(c-1)\choose e_1-e_{p+1}-c}\cdot&&\\\nonumber
[{e_1-1-e_{p+1}\choose e_{t+1}-e_{p+1}}{e_{1}-1-e_{t+1}\choose e_{t}-1-e_{t+1}}+{e_1-1-e_{p+1}\choose e_{t+1}-e_{p+1}}{e_{1}-1-e_{t+1}\choose e_{t}-e_{t+1}}]&&=\\\nonumber
\end{eqnarray}
\begin{eqnarray}\nonumber
\sum_c{e_{p+1}-1\choose c-1}{e_1-e_{p+1}-1\choose c-1}{n-e_{p+1}\choose e_1-e_{p+1}-1}{e_1-1-e_{p+1}\choose e_{t+1}-e_{p+1}}{e_1-e_{t+1}\choose e_{t}-e_{t+1}}&&=\\\label{eq:N4}
{e_1-2\choose e_1-e_{p+1}-1}{n-e_{p+1}\choose e_1-e_{p+1}-1}{e_1-e_{t+1}\choose e_{t}-e_{t+1}}{e_1-1-e_{p+1}\choose e_{t+1}-e_{p+1}}&&
\end{eqnarray}
In the first and third equality we have used part (3) of lemma~\ref{lemma:Combinatorics2}; in the last equality we have used the Vandermonde's identity.
\item If $T_{p+1}$ does not contain $n$ (by lemma~\ref{lemma:Combinatorics2} there are ${e_{p+1}-1\choose c-1}{n-e_{p+1}\choose c}$ choices for such sets) then either $T_1$ contains $(n+1)$ or it does not. Since $T_{p+1}$ does not contain $n$, there are ${n-c-e_{p+1}\choose e_1-1-c-e_{p+1}}$ choices of sets $T_1$ containing $(n+1)$. In this case either $T_{t}$ contains $(n+1)$ (there are ${e_1-1-e_{p+1}\choose e_{t}-1-e_{p+1}}$ choices of such sets) or it does not (there are ${e_1-1-e_{p+1}\choose e_{t}-e_{p+1}}$ choices of such sets). If $T_{t}$ contains $(n+1)$ then there are ${e_{t}-1-e_{p+1}\choose e_{t+1}-e_{p+1}}$ choices for $T_{t+1}$.  If  $T_{t}$ does not  contain $(n+1)$ then there are ${e_{t}-e_{p+1}\choose e_{t+1}-e_{p+1}}$ choices for $T_{t+1}$. 

The number of quadruples $\{T_1,T_{t},T_{t+1},T_{p+1}\}$ such that $T_{p+1}$ does not contain $n$ and $T_1$ contains $(n+1)$ is hence given by:
\begin{eqnarray}\nonumber
\sum_c{e_{p+1}-1\choose c-1}{n-e_{p+1}\choose c}{n-e_{p+1}-c\choose e_1-1-e_{p+1}-c}&&\\\nonumber[{e_1-1-e_{p+1}\choose e_{t}-1-e_{p+1}}{e_{t}-1-e_{p+1}\choose e_{t+1}-e_{p+1}}+{e_1-1-e_{p+1}\choose e_{t}-e_{p+1}}{e_{t}-e_{p+1}\choose e_{t+1}-e_{p+1}}]&&=\\\nonumber
\sum_c{e_{p+1}-1\choose c-1}{e_1-e_{p+1}-1\choose c}{n-e_{p+1}\choose e_1-e_{p+1}-1}{e_1-1-e_{p+1}\choose e_{t+1}-e_{p+1}}{e_1-e_{t+1}\choose e_{t}-e_{t+1}}&&=\\\label{eq:N4N1}
{e_1-2\choose e_1-e_{p+1}-2}{n-e_{p+1}\choose e_1-e_{p+1}-1}{e_1-1-e_{p+1}\choose e_{t+1}-e_{p+1}}{e_1-e_{t+1}\choose e_{t}-e_{t+1}}&&
\end{eqnarray}
By summing up \eqref{eq:N4} and \eqref{eq:N4N1} and by applying lemma~\ref{lemma:Combinatorics} we get
\begin{equation}\label{eq:89}
{e_1-1\choose e_{p+1}}{n-e_{t+1}\choose e_1-e_{t+1}-1}{n-e_{p+1}\choose e_{t+1}-e_{p+1}}{e_1-e_{t+1}\choose e_{t}-e_{t+1}}.
\end{equation}
If $T_1$ does not contain $(n+1)$ (there are ${n-c-e_{p+1}\choose e_1-c-e_{p+1}}$ choices of such sets) then there are ${e_1-e_{p+1}\choose e_{t}-e_{p+1}}$ choices for $T_{t}$ and ${e_{t}-e_{p+1}\choose e_{t+1}-e_{p+1}}$ choices for $T_{t+1}$. The number of quadruples $\{T_1,T_{t},T_{t+1},T_{p+1}\}$ such that $T_{p+1}$ does not contain $n$ and $T_1$ does not contain $(n+1)$ is hence given by:
\begin{eqnarray}\nonumber
\sum_c{e_{p+1}-1\choose c-1}{n-e_{p+1}\choose c}{n-e_{p+1}-c\choose e_1-e_{p+1}-c}{e_1-e_{p+1}\choose e_{t}-e_{p+1}}{e_{t}-e_{p+1}\choose e_{t+1}-e_{p+1}}&&=\\\nonumber
\sum_c{e_{p+1}-1\choose c-1}{e_1-e_{p+1}\choose c}{n-e_{p+1}\choose e_1-e_{p+1}}{e_1-e_{t+1}\choose e_{t}-e_{t+1}}{e_1-e_{p+1}\choose e_{t+1}-e_{p+1}}&&=\\\nonumber
{e_1-1\choose e_{p+1}}{n-e_{p+1}\choose e_1-e_{p+1}}{e_1-e_{t+1}\choose e_{t}-e_{t+1}}{e_1-e_{p+1}\choose e_{t+1}-e_{p+1}}&&=\\\label{eq:N4N1not}
{e_1-1\choose e_{p+1}}{n-e_{t+1}\choose e_1-e_{t+1}}{n-e_{p+1}\choose e_{t+1}-e_{p+1}}{e_1-e_{t+1}\choose e_{t}-e_{t+1}}&&
\end{eqnarray}
\end{enumerate}
By summing up \eqref{eq:89} and \eqref{eq:N4N1not} and by applying lemma~\ref{lemma:Combinatorics} we get the desired \eqref{M4}.

We now prove~\eqref{eq:ChiRegular}.
As for the case $p=1$ (see lemma~\ref{Lemma:RightKronecker}), the variety $Gr_\mathbf{e}(Reg^n_{p}(\lambda))$ equals the variety $Gr_\mathbf{e}(Reg^n_{p}(0))$ for every $\lambda\in k$. Indeed let us denote by $\circ$ and by $\ast$ respectively the action of $A=kQ_{p,1}$ on $Reg^n_{p}(\lambda)$ and on $Reg^n_{p}(0)$. We consider the automorphism $\psi$ of the path algebra $kQ_{p,1}$ which sends $\varepsilon_0$ to $\lambda\pi+\varepsilon_0$ where $\pi:=\varepsilon_1\circ\cdots\circ\varepsilon_p$ is the longest path of $Q_{p,1}$. As vector spaces, $Reg^n_{p}(0)$ and $Reg^n_{p}(\lambda)$ are isomorphic. Then for every $\pi$ in $A$ and every $m$ in $Reg^n_{p}(0)$, $\psi(\pi)\ast m=\pi\circ m$.  This proves that they are right--equivalent in the sense of \cite{DWZ}.
$$
\xymatrix@R=7pt@C=8pt{
                                                 &k^{n}\ar_{=}[d]&  \ar_=[l]\cdots  &k^{n}\ar_=[l]  \\
Reg^n_p(\lambda):=            &k^{n}        &                & k^n\ar^{J_n(\lambda)}[ll]\ar_{=}[u]
}
$$
We thus assume that $\lambda=0$. In this case the representation $Reg^n_{p}(0)$ is an orientable string module by lemma~\ref{Lemma:Ap1Orientable} and we can therefore apply theorem~\ref{MainThmIntro}. The Euler--Poincar\'e characteristic of $Gr_{\mathbf{e}}(Reg^n_{p}(0))$ is hence the number of $(p+1)$--tuples $\{T_{1},\cdots,T_{p+1}\}$ of subsets $T_{i}\subset [1,n]$ of cardinality $|T_{i}|=e_i$ such that $T_{i+1}\subset T_{i}$ for $i=1,\cdots, p$ and $J_n(0)(T_{p+1})\subset T_{1}$ where $J_n(0):[1,n]\rightarrow [1,n]\cup \{0\}$ is the map which sends $k$ to $k-1$, $k\in[1,n]$. The choice of the couple $\{T_{1}, T_{p+1}\}$  determines the choice of ${e_1-e_{p+1}\choose e_2-e_{p+1}}$ choices for $T_{2}$. For every such choice there are ${e_2-e_{p+1}\choose e_3-e_{p+1}}$ choices for $T_3$, and so on. For every choice of $\{T_{1},T_{p+1}\}$ there are hence $\chi_\mathbf{e}([1,p+1])=\prod_{k=1}^{p-1}{e_k-e_{p+1}\choose e_{k+1}-e_{p+1}}$ choices for $\{T_{2}, \cdots, T_{p}\}$.  The number of couples $\{T_{1}, T_{p+1}\}$ equals $\chi_{(e_1,e_{p+1})}(Reg_1^n(0))$.
It remains to prove that:
$$
\chi_{(e_1,e_2)}(
Reg_{1}^n(0))={e_1\choose e_2}{n-e_{2}\choose e_1-e_2}
$$ 
which has already been noticed in proposition~\ref{prop:EulerKronecker}.

\subsection*{Acknowledgements} I thank Professor B.~Keller for his kind hospitality and for  useful discussions on this topic during my stay in Paris.  I thank  Professor J.~Schr\"oer for many conversations about orientable string modules and for his support during my stay in Bonn. I thank Professor A.~Zelevinsky for his advice on the structure of this paper. The original publication has been accepted at the Journal Of Algebraic Combinatorics and it is available at www.springerlink.com.

\bibliographystyle{plain}
\bibliography{BibliografiaGiovanni}

\end{document}